\documentclass[11pt,a4paper]{article}

\usepackage{amsmath}
\usepackage{amsthm}
\usepackage{amsfonts}
\usepackage{amssymb}
\usepackage{amsmath,amscd}

\newcommand{\C}{\mathbb{C}}

\newcommand{\Z}{\mathbb{Z}}

\newcommand{\PP}{\mathbb{P}}
\newcommand{\RR}{\mathbb{R^{+}}}
\newcommand{\Bl}{\operatorname{Bl}}

\newcommand{\NE}{\operatorname{NE}}

\newcommand{\cont}{\operatorname{cont}}
\newcommand{\Exc}{\operatorname{Exc}}

\newcommand{\OO}{{\cal O}}

\newtheorem{theo}{Theorem}
\newtheorem{prop}{Proposition}
\newtheorem{lem}{Lemma}

\def\finpreuve
{\hskip 3pt \vrule height6pt width6pt depth 0pt}

\title
{A remark on Fano 4-folds having (3,1)-type extremal contractions}
\author{Toru Tsukioka}

\begin{document}

\maketitle

\begin{abstract}
Let $X$ be the blow-up of a smooth projective 4-fold 
$Y$ along a smooth curve $C$ and let $E$ be the exceptional divisor.
Assume that $X$ is a Fano manifold 
and has an elementary extremal contraction $\varphi:X\to Z$ of (3,1)-type 
(i.e. the exceptional locus of $\varphi$ is a divisor and its image is a curve) 
such that $E$ is $\varphi$-ample. 
We show that if the exceptional divisor of $\varphi$ is smooth, 
then $Y$ is isomorphic to $\PP^{4}$ and $C$ is an elliptic curve  
of degree 4 in $\PP^{4}$.
\end{abstract}

\section{Introduction}

As an application of the extremal contraction theory, 
S.\ Mori and S.\ Mukai 
classified smooth Fano 3-folds with Picard number greater than or 
equal to 2 (\cite{MM}). We observe that many of examples in the Mori-Mukai's 
list are obtained by blowing up other smooth projective 3-folds. 
In fact, 78 types among 88 types of smooth Fano 3-folds 
with $\rho\geq 2$ have $E_{1}$-type or $E_{2}$-type extremal 
contractions. 
In \cite{BCW} the authors classified smooth Fano varieties 
(defined over $\C$) obtained 
by blowing-up a smooth point, in any dimension. A next step 
is to consider the following problem:

{\it Problem.} Let $Y$ be a smooth projective variety. 
Let $\pi:X\to Y$ be the blow-up along a smooth
curve $C$. Classify pairs $(Y,C)$ such that $X$ is Fano.

Remark that 
for the toric case, the classification is done in any dimension by \cite{S}.
 
By the Cone and Contraction Theorem, we can 
take an extremal contraction $\varphi:X\to Z$ to normal projective variety 
such that the exceptional divisor $E$ of $\pi$ is $\varphi$-ample 
(see Lemma \ref{kisolem} below). 
It is easy to show that any fiber of $\varphi$ is at most of dimension 2. 
The author studied the case where $\varphi$ 
is a del Pezzo surface fibration and gave a complete classification (\cite{T2}).
In higher dimensions, it seems difficult to classify the case where $\varphi$ 
is birational. However, in dimension 4, 
there are several results on the birational 
extremal contractions, which may be applied to solve our problem.

In this paper, we investigate the case where $\varphi$ is of 
$(3,1)$-type contraction. 
Recall that in general, an extremal contraction 
$\varphi:X \to Z$ is said to be
$(a,b)$-type, if $\dim($Exc$(\varphi))=a$ and 
$\dim(\varphi($Exc$(\varphi)))=b$.
So, a $(3,1)$-type contraction for a 4-fold is a birational contraction 
which contracts a divisor $F$ to a curve $B$. 
The extremal contractions of $(3,1)$-type for 
smooth 4-folds are completely classified by \cite{Tk}. 
In particular, it is shown that the exceptional divisor $F$ is normal 
and $B$ is smooth. Moreover, 
$\varphi_{|F}:F\to B$ is either a $\PP^{2}$-bundle or a $Q_{2}$-bundle 
(see \cite{Tk} Main Theorem).  
\footnote{Remark that our $F$ and $B$ correspond to $E$ and $C$ in \cite{Tk}.}

In section \ref{example}, we first give an example. Let $C\subset \PP^{4}$ be 
a smooth complete intersection of one hyperplane and two hyperquadrics. 
Then, we see that $X=\Bl_{C}(\PP^{4})$ has a (3,1)-type extremal contraction 
to a complete intersection of two hyperquadrics (singular along a line) 
in $\PP^{6}$. The section \ref{proof} is devoted to show that 
this is the only example if we assume that Exc$(\varphi)$ is smooth. 
More precisely, we prove the following:

\begin{theo}\label{theorem}
Let $\pi:X\to Y$ be the blow-up of a smooth projective 4-fold 
$Y$ defined over $\C$, along a smooth curve $C$.
Assume that $X$ is a Fano manifold and 
has an elementary extremal contraction $\varphi:X\to Z$ of (3,1)-type 
such that the exceptional divisor $E$ of $\pi$ is $\varphi$-ample.
Let $F$ be the exceptional divisor of $\varphi$. 
If $F$ is smooth, then $Y$ is isomorphic to $\PP^{4}$ 
and $C$ is a smooth complete intersection of a hyperplane and 
two hyperquadrics.
\end{theo}

We will use the following lemma, which is essentially 
the same as in \cite{BCW}(Lemme 2.1). For reader's convinience, 
we include here the statement with its proof.
\begin{lem}\label{kisolem}
Let $X$ be a Fano manifold and let $E$ be a non-zero 
effective divisor on $X$. Then there exists an extremal ray 
$\RR[f]\subset\overline{\NE}(X)$ such that $E\cdot f>0$.
\end{lem}

{\it Proof.} Since $X$ is projective, we can take a 
curve $\Gamma$ on $X$ such that $E\cdot \Gamma>0$. 
By the Cone Theorem, there exist positive real numbers $a_{i}$, 
and extremal rational curves $f_{i}$ 
such that
$\Gamma\equiv \sum a_{i}f_{i}$ (finite sum). 
Hence
$$
0<E\cdot \Gamma=\sum a_{i}(E\cdot f_{i}). 
$$
This implies that one of extremal rational curves satisfies 
$E\cdot f_{i}>0$. \finpreuve

\

Throughout this paper, we shall assume that the base field is 
the complex numbers.
For a Cartier divisor $D$ and a 1-cycle $\alpha$ on a variety $X$, 
we denote the intersection number by $D\cdot \alpha$, 
but we also write $(D\cdot \alpha)_{X}$ when we need to   
clarify the variety in which the intersection number is taken.

\section{An example}\label{example}

We give an example of a smooth Fano 4-fold $X$ obtained by blowing up 
along a curve 
such that $X$ has another $(3,1)$-type extremal contraction.

{\bf Example} Let $C\subset \PP^{4}$ be a smooth complete intersection 
of a hyperplane and two hyperquadrics, $\pi:X\to \PP^{4}$ 
the blow-up along $C$, and $E$ the exceptional divisor.
Let $F$ be the strict transform of the hyperplane containing $C$. 
Remark that $F\simeq \Bl_{C}(\PP^{3})$ is a $Q_{2}$-bundle over 
$\PP^{1}$.
Let $e$ be a line in a fiber of the $\PP^{2}$-bundle 
$\pi_{|E}:E\to C$, and let $f$ be the strict transform 
of a line in $\PP^{4}$ intersecting $C$ at two points.
Then we have 
$$\overline{\NE}(X)=\RR[e]+\RR[f].$$

The extremal contraction 
associated to the ray $\RR[e]$ is of course the blow-up 
$\pi:X\to \PP^{4}$. 
Let $L:=\pi^{*}\OO_{\PP^{4}}(1)$. The linear system $|2L-E|$ is base-point-free and defines 
the extremal contraction $\varphi:X\to Z$ of the ray $\RR[f]$. 
Indeed, we have $(2L-E)\cdot f=0$. Note that $B:=\varphi(F)$ is 
isomorphic to $\PP^{1}$ and $\varphi_{|F}:F\to B$ is a $Q_{2}$-bundle. 
Thus, $\varphi$ is a (3,1)-type extremal contraction whose exceptional divisor is $F$.
More precisely, the image $Z$ is a complete intersection of two hyperquadrics in $\PP^{6}$, 
singular along $B\simeq\PP^{1}$.
To see this, we calculate $h^{0}(X,\OO_{X}(2L-E))$ and $(2L-E)^{4}$.

Consider the exact sequence:
$$
0\to \OO_{X}(2L-E)\to\OO_{X}(2L)\to \OO_{E}(2L)\to 0.
$$
Remark that $A:=-K_{X}+(2L-E)=(5L-2E)+(2L-E)=7L-3E$ is ample 
by Kleiman's criterion, 
because $A\cdot e=3>0$ and $A\cdot f=1>0$. Therefore, by the Kodaira vanishing, $H^{1}(X,\OO_{X}(2L-E))=0$.
On the other hand, we get $h^{0}(X,\OO_{X}(2L))=h^{0}(\PP^{4},\OO_{\PP^{4}}(2))=15$. Since $\OO_{E}(2L)\simeq(\pi_{|E})^{*}\OO_{C}(2)$, we have
$h^{0}(E,\OO_{E}(2L))=h^{0}(C,\OO_{C}(2))=\deg(\OO_{C}(2))=8$ 
(recall that $\pi_{|E}$ is a $\PP^{2}$-bundle and $g(C)=1$). 
Hence, 
$$
h^{0}(X,\OO_{X}(2L-E))=h^{0}(X,\OO_{X}(2L))-h^{0}(E,\OO_{E}(2L))=7
$$ 
and $|2L-E|$ defines a morphism $\varphi:X\to \PP^{6}$. Now we determine the 
image of $X$.
Note that we have $L^{2}\cdot E\equiv 0$, $L\cdot E^{3}=\deg C=4$, and 
$E^{4}=\deg N_{C/\PP^{4}}=20$. Thus, 
$$
(2L-E)^{4}=(2L)^{4}-8L\cdot E^{3}+E^{4}=4.
$$
Consider the exact sequence
$$
0\to \OO_{\PP^{6}}(2)\otimes I_{Z}\to\OO_{\PP^{6}}(2)\to\OO_{Z}(2)\to 0.
$$
Since $h^{0}(Z,\OO_{Z}(2))=h^{0}(X,\OO_{X}(4L-2E))=26$, we obtain
$$
h^{0}(\PP^{6},\OO_{\PP^{6}}(2)\otimes I_{Z})
\geq h^{0}(\PP^{6},\OO_{\PP^{6}}(2))-h^{0}(Z,\OO_{Z}(2))=28-26=2.
$$
It follows that there exist two linearly independent 
hyperquadrics in $\PP^{6}$ containing $Z$. 
Since $\deg Z=(2H-E)^{4}=4$, $Z$ is a complete intersection of two hyperquadrics.

\section{Proof of Theorem \ref{theorem}}\label{proof}

Denote by $e$ a line in a fiber of the 
$\PP^{2}$-bundle $\pi_{|E}:E\to C$. 
The key to prove Theorem \ref{theorem} is the following:

\begin{lem}\label{lemma}
We have $F\cdot e=1$.
\end{lem}

{\it Proof.} We denote by $(e)$ the corresponding point 
in Hilb$(X)$. Let $T$ be the reduced part of 
the irreducible component of 
Hilb$(X)$ containing $(e)$. Note that 
$T$ is a $\PP^{2}$-bundle over $C$ whose fiber $T_{c}$ ($c\in C$) 
parametrizes lines in $E_{c}:=\pi^{-1}(c)\simeq \PP^{2}$. In particular, 
$T$ is smooth and of dimension 3.

{\it Step 1.} For all $(e)\in T$ such that $e\not\subset F$, 
we have $\sharp(F\cap e)=1$. 
\footnote{We mean by $\sharp(F\cap e)$ the number of points on $F\cap e$ 
without multiplicity.}
Assume the contrary, i.e. there exists $(e_{0})\in T$ 
such that $e_{0}\not\subset F$ and $\sharp(F\cap e_{0})\geq 2$. 
Remark that $\varphi(e_{0})\neq B$. Let $x_{i}$ $(i=1,2)$ be 
two distinct points in $F\cap e_{0}$ and let $b_{i}:=\varphi(x_{i})$.
Consider the incidence graph:
$$
\begin{CD}
V   @>p >> X \\
@VqVV \\
T
\end{CD}
$$
We define $V_{i}:=p^{-1}(E\cap \varphi^{-1}(b_{i}))$ and 
$T_{i}:=q(V_{i})$ for $i=1,2$. 
Note that $\dim V_{i}=\dim (E\cap \varphi^{-1}(b_{i}))+1$  
because $p$ is a $\PP^{1}$-bundle. We observe that 
$q_{|V_{i}}$ is a finite map. Indeed, if not, 
there exists $t\in T_{i}$ such that $q^{-1}(t)\subset V_{i}$. 
Then $e_{t}:=p(q^{-1}(t))$ is contracted by $\varphi$. 
This contradicts to our assumption that $E$ is $\varphi$-ample. 
It follows that $\dim T_{i}=\dim V_{i}=2\ (i=1,2)$. 
Note also that $(e_{0})\in T_{1}\cap T_{2}$. Now, we have
$$
\dim(T_{1}\cap T_{2})\geq \dim T_{1}+\dim T_{2}-\dim T=2+2-3=1.
$$ 
So, we can take an irreducible curve $A\subset T_{1}\cap T_{2}$ 
passing through $(e_{0})$. Then $q^{-1}(A)$ is a ruled surface 
having two exceptional curves $V_{i}\cap S\ (i=1,2)$, a contradiction.

{\it Step 2.} Consider $M:=(F\cap E)_{red}$. 
By Step 1, we see that for each $c\in C$, $e_{c}:=(F\cap E_{c})_{red}$ is a line in 
$E_{c}\simeq \PP^{2}$. 
So, $\pi_{|M}:M\to C$ is a $\PP^{1}$-bundle. In particular $M$ is irreducible. 
We can write 
$E_{|F}=mM$ with $m\in \Z^{+}$. We have
$$
(mM\cdot e_{c})_{F}=(E_{|F}\cdot e_{c})_{F}
=(E\cdot e_{c})_{X}=-1
$$  
 
By assumption, $F$ is smooth. So, $M\subset F$ is a Cartier divisor and 
$(M\cdot e_{c})_{F}$ is integer. It follows that $m=1$, i.e. the 
intersection $F\cap E$ is transversal. We conclude that 
$F\cdot e=\sharp(F\cap e)=1$. \finpreuve

\

{\it Proof of Theorem \ref{theorem}.} \ 
By the proof of Lemma \ref{lemma}, 
$\pi_{|M}:M\to C$ is a $\PP^{1}$-bundle and
$(M\cdot e_{c})_{F}=-1$.
So, $\pi_{|F}:F\to F':=\pi(F)$ is 
the blow-up along $C$, and $F'$ is smooth. 
On the other hand, by \cite{Tk}, $\varphi_{|F}:F\to B$ is either 
a $\PP^{2}$-bundle or a $Q_{2}$-bundle. 
Therefore $F$ is a Fano 3-fold with $\rho(F)=2$. By assumption, $F$ is smooth. 
So, by the Mori-Mukai's list, 
the pair $(F',C)$ is one of the following: 
\begin{itemize}
\item[(1)] $F'\simeq \PP^{3}$ and $C$ is a line;
\item[(2)] $F'$ is a hyperquadric $Q_{3}\subset\PP^{4}$ and $C=H\cap H'$ 
with $H$, $H'\in |\OO_{Q_{3}}(1)|$;
\item[(3)] $F'\simeq \PP^{3}$ and $C=Q\cap Q'$ with 
$Q$, $Q'\in |\OO_{\PP^{3}}(2)|$.
\end{itemize}
In the case (3), $C$ is an elliptic curve. So, 
$Y$ is a Fano manifold by \cite{Wis} (Proposition 3.5). 
In the cases (1) and (2), we have $N_{C/F'}\simeq \OO_{C}(1)^{\oplus 2}$.
Since there exists an inclusion of normal bundles $N_{C/F'}\subset N_{C/Y}$, 
$N_{C/Y}$ cannot be isomorphic to 
$\OO_{\PP^{1}}(-1)^{\oplus 3}$. So, $Y$ is a Fano manifold 
again by \cite{Wis}. 

Now, by Lemma \ref{kisolem}, we can take an extremal ray $\RR[m]$ 
such that $F'\cdot m>0$. 
Then, by Proposition \ref{proposition} below, we have $\rho(Y)=1$.
In paticular $F'$ is ample. 
Let $f$ be a minimal rational curve of 
the extremal contraction $\varphi$. 
We obtain the following table of intersection numbers 
(due to \cite{Tk} and \cite{MM}):
$$
\begin{tabular}{c|c|c}
\hline
case& $F\cdot f$ & $E\cdot f$ \\ \hline \hline
(1)& $-1$ or $-2$& $1$\\
(2)& $-1$        & $1$\\
(3)&$-1$         & $2$\\
\hline
\end{tabular} 
$$

Let $f':=\pi_{*}f$. Note that 
$F'\cdot f'=(\pi^{*}F')\cdot f=(F+E)\cdot f$. 
In the cases (1) and (2), we have $F'\cdot f'\leq 0$, a contradiction 
because $F'$ is ample. 
So, only the case (3) (in which we have $F'\cdot f'=1$) is possible, and $(Y,F')\simeq (\PP^{4},\OO_{\PP^{4}}(1))$. 
Consequently, $C$ is the complete intersection 
$F'\cap Q\cap Q'$ with $F'\in |\OO_{\PP^{4}}(1)|$ and 
$Q, Q'\in |\OO_{\PP^{4}}(2)|$. \finpreuve

\

It remains to prove the following:

\begin{prop}\label{proposition}
Let $Y$ be a smooth projective variety of dimension $n\geq 4$ and 
$D$ a prime divisor on $Y$ with $\rho(D)=1$. 
Assume that there exists an extremal contraction 
$\mu:Y\to V$ of ray $\RR[m]$ with $D\cdot m>0$, 
$m$ being a minimal rational curve of the ray.
If there exsits a smooth curve $C\subset D$ such that the blow-up 
$X:=\Bl_{C}(Y)$ is a Fano manifold, then we have $\rho(Y)=1$. 
Moreover, if $D$ is isomorphic to $\PP^{n-1}$, then we have 
$(Y,D)\simeq (\PP^{n},\OO_{\PP^{n}}(1))$.
\end{prop}

{\it Proof.} We shall consider two cases:
\begin{itemize}
\item[(1)] there exists $v_{0}\in V$ such that $\dim (\mu^{-1}(v_{0})\cap D)\geq 1$;
\item[(2)] $\dim (\mu^{-1}(v)\cap D)=0$ for all $v\in V$.
\end{itemize}
In the case (1), there exsits a curve $B\subset \mu^{-1}(v_{0})\cap D$. 
So, we can write $B\equiv bm$ with $b\in \RR$. 
Since $\rho(D)=1$, any curve in $D$ is numerically equivalent to a multiple of $m$. 
Hence, $\mu(D)$ is a point. 
We also have $D\cdot B>0$. 
Now, by Proposition 4 of \cite{T2}, 
we conclude that $\rho(Y)=1$.

We show that the case (2) is impossible. 
In this case, any fiber of $\mu$ is at most of dimension 1. So, by \cite{Ando} (see 
also \cite{Wis} Theorem 1.2), $\mu$ is either, a $\PP^{1}$-bundle, 
a conic bundle, or a blow-up along  
a smooth subvariety of codimension 2 in a smooth projective variety. 
If $\mu$ is a $\PP^{1}$-bundle, take a fiber $m$ passing 
through a point on $C$.
Let $\tilde{m}$ be the strict transform by the blow-up $\pi:X\to Y$. 
For the exceptional divisor $E$, we have $E\cdot \tilde m\geq 1$, so that
$$
K_{X}\cdot \tilde{m}=K_{Y}\cdot m+(n-2)E\cdot\tilde{m}\geq -2+(n-2)=n-4\geq 0,
$$  
which is absurd because $X$ is a Fano manifold.

If $\mu$ is a conic bundle, the extremal rational curve $m$ 
is a component of a singular fiber of $\mu$. 
Let $\Delta$ be the discriminant locus and let 
$\tilde{\Delta}:=\mu^{-1}(\Delta)$. 
The assumption $D\cdot m>0$ implies 
$\tilde{\Delta}\cap D \neq \varnothing$. Since $\rho(D)=1$, the non-zero effective 
Cartier divisor $\tilde{\Delta}_{|D}$ is ample. 
Therefore, 
$$
(\tilde{\Delta}\cdot C)_{Y}=(\tilde{\Delta}_{|D}\cdot C)_{D}>0, 
$$
so that $\tilde{\Delta}\cap C\neq \varnothing$. Now, we can take 
a sigular fiber $\mu^{-1}(v_{0})$ ($v_{0}\in \Delta$) meeting $C$.
Let $m_{0}\subset\mu^{-1}(v_{0})$ be a component 
such that $m_{0}\cap C\neq \varnothing$. 
Then, we have a contradiction as in the case of 
$\PP^{1}$-bundle. 
The case of a blow-up along a centre of codimension 2, can be ruled out 
by using a same argument for the exceptional divisor of $\mu$ 
in place of $\tilde{\Delta}$.

Consequently, only the case (1) is possible, so that we have $\rho(Y)=1$. 
If $D\simeq\PP^{n-1}$, by \cite{BCW}(Lemme 4) 
we conclude that $(Y,D)\simeq (\PP^{n},\OO_{\PP^{n}}(1))$ 
\finpreuve    

\

Our assumption that $F=\Exc (\varphi)$ is smooth, 
is used in the proof of Lemma \ref{lemma} (only for {\it Step.2}) and 
in the proof of Theorem \ref{theorem} in oder to apply to $F$ the Mori-Mukai's classification 
of smooth Fano 3-folds.
So, it is natural to ask whether Theorem \ref{theorem} remains true 
without the smoothness of $F$.
Concerning to this question, it is worth seeing 
the following:

{\bf Example} (A degenerate case of the example in Section \ref{example})
We consider the union of two smooth conics $C=C_{1}\cup C_{2}\subset Y:=\PP^{4}$
obtained as complete intersection 
of a hyperplane and two hyperquadrics. 
We assume that $C_{1}$ and $C_{2}$ meet 
at two distinct points.
Let $\pi:X\to \PP^{4}$ be 
the blow-up along the ideal $I_{C_{1}\cup C_{2}}$ and $E$ the exceptional divisor.
Let $F$ be the strict transform of the hyperplane containing $C=C_{1}\cup C_{2}$. 
Then $F$ is a $Q_{2}$-bundle over 
$\PP^{1}$ having exactly two ordinary double points. Remark that $F$ 
is isomorphic to the blow-up of $\PP^{3}$ along the ideal $I_{C_{1}\cup C_{2}}$. 
Moreover, $F$ can be realized as divisor in $\PP^{1}\times \PP^{3}$ 
by the equation $sX_{2}X_{3}+t(X_{0}^{2}+X_{1}^{2}+X_{2}^{2}+X_{3}^{2})=0$, 
where $(s:t)$ (resp. $(X_{0}:X_{1}:X_{2}:X_{3})$) is the homogeneous coodinates 
of $\PP^{1}$ (resp. $\PP^{3}$). The fiber over $(1:0)$ is two planes 
$P_{i}$ $(i=1,2)$ and the two ordinary double points lie on the line $P_{1}\cap P_{2}$.
  
As in Section \ref{example}, 
we see that the linear system $|\pi^{*}\OO_{\PP^{4}}(2)-E|$ defines a (3,1)-type contraction 
$\varphi:X \to Z$ to complete intersection of two hyperquadrics in $\PP^{6}$, and 
its exceptional divisor is $F$. This gives an example of $(Y,C)$ such that 
$F=\Exc(\varphi)$ is singular. However $X$ is also singular along two rational curves 
over the two intersection points of $C_{1}\cap C_{2}$.

\section{Related results} 

Let $X$ be a Fano manifold and let $\iota_{X}$ 
be its pseudo-index, i.e. the minimum of the anti-canonical degrees 
$(-K_{X}\cdot C)$ for rational curves $C$ on $X$. 
In \cite{BCDD}, the authors discuss the inequality ("generalized Mukai conjecture"): 
$$
\rho(X)(\iota_{X}-1)\leq \dim X
$$
and prove it in dimension 4.
The essential part is to show that if $\iota_{X}=2$, then $\rho(X)\leq 4$.
Concerning to this, we have the following:  

\begin{prop}\label{intersection}
Let $\pi:X\to Y$ be the blow-up of a smooth projective variety 
$Y$ of dimension $n\geq 4$ 
along a smooth curve $C$ and let $E$ be the exceptional divisor.
Assume that $X$ is a Fano manifold 
and there is another blow-up $\varphi:X\to Z$ (different from 
$\pi$) along a smooth 
curve $B$. Let $F$ be the exceptional divisor of $\varphi$. Then, 
we have $E\cap F=\varnothing$.

\end{prop} 
{\it Proof.} Assume $E\cap F\neq \varnothing$. Take $a\in C$ and 
$b\in B$ such that $E_{a}\cap F_{b}\neq\varnothing$. Then we obtain
$\dim (E_{a}\cap F_{b})\geq \dim E_{a}+\dim F_{b}-\dim X=n-4$. 
So, if $n\geq 5$, there is a curve contained in $E_{a}\cap F_{b}$ 
and then contracted by both $\pi$ and $\varphi$. This is absurd 
because we assume $\pi\neq \varphi$. Therefore, we have $n=4$. 
By (the proof of) Theorem \ref{theorem}, 
$\varphi_{|F}:F\to B$ cannot be a $\PP^{2}$-bundle. So, the case 
$E\cap F\neq \varnothing$ is impossible.
\finpreuve

\

We are now able to give a simple proof of a result in \cite{BCDD}.

\begin{theo}[see \cite{BCDD} Th\'eor\`eme 3.9]
Let $X$ be a Fano manifold of dimension $\geq 4$ whose 
birational contractions are all blow-ups along 
smooth curves in smooth projective varieties. Assume 
that $X$ has at least one birational contraction. Then, we have 
$\rho(X)\leq 3$.
\end{theo}
{\it Proof.} 
Let $E$ be an exceptional divisor on $X$.
By Lemma \ref{kisolem}, we can take an extremal ray 
$\RR[f]\subset \overline{\NE}(X)$ such that 
$E\cdot f>0$ . Then, by assumption and by 
Proposition \ref{intersection} above, the associated contraction 
$\mu:=\cont_{\RR[f]}:X\to Z$ is of fiber type. 
So, there is a surjection $\mu_{|E}:E\to Z$.
Hence, we have $\rho(Z)\leq\rho(E)=2$. 
Consequently, $\rho(X)=\rho(Z)+1\leq 3$.
\finpreuve

\

\flushright{
Department of Mathematics \ \ \ \ \ \ \\
Tokyo Institute of Technology \ \ \ \ \\
2-12-1 Oh-okayama, Meguro-ku, \ \\
Tokyo 152-8551, JAPAN \ \ \ \ \ \ \ \ \ \ 

\

email: tsukiokatoru@yahoo.co.jp

}

\end{document}